\documentclass[a4paper,12pt]{article} 
\usepackage[parfill]{parskip}
\usepackage[titletoc,toc,title]{appendix}
\usepackage[margin=2cm]{geometry}
\usepackage{ mathrsfs }
\usepackage{amsthm}
\usepackage{amsmath}
\usepackage{amssymb}
\usepackage{graphicx}
\usepackage{mathtools}
\usepackage{subcaption}
\usepackage{enumerate}
\usepackage{enumitem}
\usepackage{xcolor}
\usepackage{tikz}
\usetikzlibrary{calc,intersections,through,backgrounds}
\usepackage{array}
\usepackage{tikz-cd}

\usepackage[bookmarks=true]{hyperref}
\hypersetup{
	colorlinks,
	linkcolor={red!70!black},
	citecolor={red!90!black},
	urlcolor={green!80!black}
	linktoc=all, 
}

\makeatletter
\newtheorem*{rep@theorem}{\rep@title}
\newcommand{\newreptheorem}[2]{%
	\newenvironment{rep#1}[1]{%
		\def\rep@title{#2 \ref{##1}}%
		\begin{rep@theorem}}%
		{\end{rep@theorem}}}
\makeatother
\newtheoremstyle{stylename}
{15pt} 
{15pt} 
{\itshape} 
{} 
{\bfseries} 
{.} 
{.5em} 
{} 
\theoremstyle{stylename}
\newtheorem{theorem}{Theorem}[section]
\newreptheorem{theorem}{Theorem}
\newtheorem{lemma}[theorem]{Lemma}
\newtheorem{corollary}[theorem]{Corollary}

\newtheoremstyle{exampstyle}
{15pt} 
{15pt} 
{} 
{} 
{\bfseries} 
{.} 
{.5em} 
{} 
\theoremstyle{exampstyle}

\newtheorem{remark}[theorem]{Remark}

\usepackage{titlesec}

\titleformat{\chapter}[display]
  {\normalfont\sffamily\huge\bfseries\color{black}}
  {\chaptertitlename\ \thechapter}{20pt}{\Huge}

\titleformat{\section}
  {\Large\bfseries\color{black}}
  {\thesection}{1em}{}

\begin{document}


\title {Three-dimensional  connected groups of automorphisms of \\ toroidal circle planes}
\author{Brendan Creutz, Duy Ho, G\"{u}nter F. Steinke}
 \maketitle
	
	\begin{abstract}
		We contribute to the classification of toroidal circle planes and flat Minkowski planes possessing three-dimensional connected groups of automorphisms. When such a group is an  almost simple Lie group, we show that it is isomorphic to $\text{PSL}(2,\mathbb{R})$. 
		Using this result, we describe a framework for the 	full classification based on the action of the group on the point set.
		 
\end{abstract}


\section{Introduction}

In  incidence geometry, Minkowski planes are one of the three types of Benz planes, the other two are M\"{o}bius (inversive) and Laguerre planes.  
There are different ways to consider Minkowski planes. They appear in the context of  finite geometries, with  applications in coding theory, cf. \cite{karzel1990}. 
In the work of Benz \cite{benz1973}, Minkowski planes were studied  from an algebraic point of view. Perhaps this and the influence of Salzmanns's work on $\mathbb{R}^2$-planes (cf. \cite{salzmann1967b}) leads to the development of Minkowski planes in a topological setting. 

Under suitable topological assumptions, a Minkowski plane  defined over the reals is called a flat Minkowski plane (cf. Subsection 2.1).
In contrast to the scarcity of known models defined over other fields, there are many examples of flat Minkowski planes, cf. \cite{gunter2001} or  \cite{duythesis} and references therein. 
The classical model of a flat Minkowski plane is the geometry of plane sections of the standard nondegenerate ruled quadric in real 3-dimensional projective space. 
An algebraic description of the classical flat Minkowski plane can be found in \cite[Subsection 4.1.5]{gunter2001}. 

Schenkel \cite{schenkel1980} showed that the automorphism group of a flat Minkowski plane is a Lie group with respect to the compact-open topology and has dimension at most 6. 
One way to investigate flat Minkowski planes  is to describe all possible planes with automorphism groups of a given dimension $n$. For convenience, we call this number $n$ the \textit{group dimension} of the plane.  It is known that planes of group dimension at least 5 are isomorphic to the classical real Minkowski plane, which has group dimension 6. Also in \cite{schenkel1980}, Schenkel  determined planes of group dimension 4 or with 3-dimensional kernels. In case of flat Minkowski planes of group dimension 3, despite many examples, a full classification for this dimension is not complete.

Toroidal circle planes are  a generalisation of flat Minkowski planes in the sense that they are  required to satisfy  all but one geometric axiom for flat Minkowski planes. For precise definitions of these incidence structures and related terminology, we refer the reader to Subsection 2.1.  So far, the only examples of  toroidal circle planes that are not flat Minkowski planes are those constructed by Polster \cite{polster1998b}.
 In \cite{brendan2017a}, we showed that the automorphism group of a toroidal circle plane is also a Lie group with  dimension at most 6. This allows us to extend the ongoing  classification of flat Minkowski planes to the more general toroidal circle planes. 

Also in the same paper above, we showed that toroidal circle planes of group dimension at least 4 or with 3-dimensional kernels are precisely the flat Minkowski planes that Schenkel described. As it can be applied for any group dimension, the machinery used by Schenkel plays an important role in our framework. The action of a  group of automorphisms $G$ is dictated by its action on each factor $\mathbb{S}^1$ of the torus. Groups acting on 1-manifolds are known by Brouwer's Theorem (cf. Theorem \ref{brouwer}), and when the group $G$ is large enough, its structure can be fully described. This is how Schenkel's result (for flat Minkowski planes)  and \cite[Theorem 1.2]{brendan2017a} (for toroidal circle planes) were obtained.

Concerning the next open case of toroidal circle planes with 3-dimensional  groups of automorphisms, we face the following two objectives:

(A) Narrow down the possibilities for  a 3-dimensional connected  group of automorphisms of a toroidal circle plane.

(B) Given a group $G$ that can possibly occur as in (A), determine the existence and  characterise toroidal circle planes admitting $G$ as its group of automorphisms. 

In this paper, we present results in the direction of (A) thereby laying the groundwork to address  (B) systematically in the future.  
The  paper is organized as follows.  From Brouwer's Theorem, we first obtain an initial list of possibilities for a 3-dimensional group. Next, in the special case when such a group is an almost simple Lie group, we show that it is  isomorphic to $\text{PSL}(2,\mathbb{R})$. We then   determine how $\text{PSL}(2,\mathbb{R})$ acts as a group of automorphisms. 
This result on almost simple Lie groups is independent from the group dimension,  and we state it as our first main theorem.

\begin{theorem} \label{almostsimplePSL}
		Let $S$ be an almost simple connected group of automorphisms of a toroidal circle plane $\mathbb{T}$. Then $S \cong \text{\normalfont PSL}(2,\mathbb{R})$. Furthermore, exactly one of the following occurs.
			\begin{enumerate}
				\item $S$ fixes either every $(+)$-parallel class or  every $(-)$-parallel class. In both cases, $\mathbb{T}$ is isomorphic to a half-classical Minkowski plane $\mathcal{M}(f,id)$ (described below), where $f$ is an orientation-preserving homeomorphism of $\mathbb{S}^1$. 
				
				\item  $S$ acts diagonally on the point set. The diagonal $D$ is a circle, and $S$ fixes $D$. 
			\end{enumerate}
\end{theorem}


Half-classical Minkowski planes can be described briefly as follows. Let $f$ and $g$ be two orientation-preserving homeomorphisms of $\mathbb{S}^1$.
The circle set $\mathcal{C}(f,g)$ of a \textit{half-classical  Minkowski plane $\mathcal{M}(f,g)$} consists of sets of the form 
$$
\{ (x,\gamma(x)) \mid x \in \mathbb{S}^1 \},
$$
where $\gamma \in \text{\normalfont PSL}(2,\mathbb{R}) \cup g^{-1}(\text{\normalfont PGL}(2,\mathbb{R}) \backslash \text{\normalfont PSL}(2,\mathbb{R}))f$. For a reference, cf. \cite[p. 239]{gunter2001}.

With the aid of Theorem \ref{almostsimplePSL} for almost simple groups and previous work in the literature for other cases,  we determine all possible geometric configurations a 3-dimensional group can fix. 
We denote the connected component of the affine group $\text{AGL}_1(\mathbb{R})$ by
$$
\text{\normalfont L}_2 = \{x \mapsto ax+b \mid a, b \in \mathbb{R}, a >0\}.
$$
We also define the following subgroups of $\text{AGL}_2(\mathbb{R})$:
$$
\Phi_\infty \coloneqq \{(x,y) \mapsto (x+b,ay+c) \mid a,b,c \in \mathbb{R}, a>0\},
$$
and
 $$
 \Phi_d \coloneqq \{(x,y) \mapsto (ax+b,a^dy+c) \mid a,b,c \in \mathbb{R}, a>0\},
 $$
 for $d \in \mathbb{R}$. In the second half of the paper, we prove the following. 
 
\begin{theorem} \label{3dfull} Let $\Sigma$ be a 3-dimensional connected group of automorphisms of a toroidal circle plane $\mathbb{T}=(\mathcal{P}, \mathcal{C}, \mathcal{G}^+, \mathcal{G}^-)$. Let $\Delta^\pm$ be the kernel of the induced action of $\Sigma$ on $\mathcal{G}^\pm$. Then exactly one of the following occurs.
	\begin{enumerate}
		
			\item   $\Sigma$ fixes no parallel classes but fixes and acts transitively on exactly one circle. In this case $\Sigma \cong \text{\normalfont \text{\normalfont PSL}}(2,\mathbb{R})$ and it acts diagonally on the torus, under suitable coordinates. The diagonal is the fixed circle under this action. 
			
				\item  $\Sigma$ fixes no points but fixes and acts transitively on either every $(+)$-parallel class  or every $(-)$-parallel class. In this case $\Sigma \cong \text{\normalfont \text{\normalfont PSL}}(2,\mathbb{R})$ and $\mathbb{T}$ is isomorphic to a half-classical Minkowski plane $\mathcal{M}(f,id)$, where $f$ is an orientation-preserving homeomorphism of $\mathbb{S}^1$.

			\item  $\Sigma$ fixes no points but fixes and acts transitively on exactly one parallel class $\pi$. 
			In this case  $\Sigma \cong  \text{\normalfont L}_2 \times \text{\normalfont SO}(2,\mathbb{R})$. 
			Assume $\pi$ is a $(+)$-parallel class.  
			Then the factor group $\Sigma /\Delta^-$  
			is isomorphic   
			and  acts equivalently to $\text{\normalfont SO}(2,\mathbb{R})$
			on $\mathcal{G}^-$. Also, 	the factor group  $\Sigma /\Delta^+ $
			is isomorphic   
			and  acts equivalently to  $\text{\normalfont L}_2 $
			on   $\mathcal{G}^+\backslash\{\pi\}$. 
			 The case when $\pi$ is a $(-)$-parallel class is analogous. 
			
					\item $\Sigma$ fixes exactly two parallel points.  The coordinates may be chosen such that the two  fixed points are $(\infty,\infty)$ and $(0,\infty)$. In this case $\Sigma \cong \Phi_d$, for some $d \le 0$,  and the  action of $\Sigma$ is described by the maps
			$$
			\{ (x,y) \mapsto (ax,by+c) \mid a,b>0,c \in \mathbb{R} \},
			$$
			when $\Sigma \cong \Phi_0$, and  
			$$
			\{ (x,y) \mapsto (a \text{ } \text{sgn}(x) \cdot |x|^b ,b^dy+c) \mid a,b>0,c \in \mathbb{R} \},
			$$
			when $\Sigma \cong \Phi_d$, when $d<0$. 
			
			\item$\Sigma$ fixes exactly one point. The coordinates may be chosen such that the fixed point is $p=(\infty,\infty)$.  In this case the derived plane $\mathbb{T}_p$ is Desarguesian and $\Sigma \cong \Phi_d$, for some $d \in \mathbb{R} \cup \{ \infty\}$. The action of  $\Sigma$ is described by the standard action of $\Phi_d$ on $\mathbb{R}^2$.

	\end{enumerate} 
\end{theorem} 

 \begin{remark} \label{remarkexample} 
	
	Besides Case 2 which is fully determined, there are examples of flat Minkowski planes of group dimension 3 for  Cases 1 and 5. Two families of flat Minkowski planes admitting 3-dimensional groups fixing no points but fixing and acting transitively on a circle were constructed by Steinke \cite{gunter2004} and \cite{gunter2017}. An Artzy-Groh plane  $\mathcal{M}_{AG}(f,g)$  (cf. \cite{artzy1986})   admits the group $\Phi_1$ with its standard action.  
	It is currently unknown if there are toroidal circle planes satisfying the conditions in Cases 3 and 4.
\end{remark}
 

\section{Preliminaries}

\subsection{ Toroidal circle planes, flat Minkowski planes and derived planes}

A \textit{toroidal circle plane} is a geometry $\mathbb{T}=(\mathcal{P}, \mathcal{C}, \mathcal{G}^+, \mathcal{G}^-)$, whose 
\begin{enumerate}[label=$ $]
	\item point set $\mathcal{P}$ is  the torus $\mathbb{S}^1 \times \mathbb{S}^1$, 
	\item circles (elements of $\mathcal{C}$) are graphs of homeomorphisms of $\mathbb{S}^1 $,  
	\item $(+)$-parallel classes (elements of $\mathcal{G}^+$) are the verticals $\{ x_0 \}  \times \mathbb{S}^1$,  
	\item $(-)$-parallel classes (elements of $\mathcal{G}^-$) are the horizontals $\mathbb{S}^1 \times \{ y_0 \}$,  
\end{enumerate} 
where $x_0, y_0 \in \mathbb{S}^1$. 

We denote the $(\pm)$-parallel class containing a point $p$ by $[p]_\pm$. When two points $p,q$ are on the same $(\pm)$-parallel class, we say they are \textit{$(\pm)$-parallel}  and denote  this by $p \parallel_{\pm} q$. Two points $p,q$ are $\textit{parallel}$ if they are  $(+)$-parallel or $(-)$-parallel, and we denote this by $p \parallel q$. 

Furthermore, a toroidal circle plane satisfies the following
\begin{enumerate}[label=]
	\item \textit{Axiom of Joining}: three pairwise non-parallel points  $p,q,r$ can be joined by a unique circle $\alpha(p,q,r)$. 
\end{enumerate}

A toroidal circle plane is called a \textit{flat Minkowski plane}  if it also satisfies the following
\begin{enumerate}[label=$ $]
	\item \textit{Axiom of Touching}: for each circle $C$ and any two nonparallel points $p,q$ with $p \in C$ and $q \not \in C$, there is exactly one circle $D$ that contains both points $p,q$ and intersects $C$ only at the point $p$.
\end{enumerate} 


	The \textit{derived plane $\mathbb{T}_p$ of $\mathbb{T}$ at the point $p$} is the incidence geometry whose point set is $\mathcal{P} \backslash ( [p]_+ \cup [p]_-)$, whose lines are all parallel classes not going through $p$ and all circles of $\mathbb{T}$ going through $p$.   
	For every point $p \in \mathcal{P}$, the derived plane $\mathbb{T}_p$ is an $\mathbb{R}^2$-plane  and even a flat  affine plane   when $\mathbb{T}$ is a flat Minkowski plane, cf. \cite[Theorem 4.2.1]{gunter2001}.

\subsection{The automorphism group}
%

An \textit{automorphism of a   toroidal circle plane $\mathbb{T}$} is a permutation of the point set $\mathcal{P}$ such that parallel classes are mapped to parallel classes and circles are mapped to circles.    With respect to composition, the set of all automorphisms of a toroidal circle plane is an abstract group, which we call the \textit{automorphism group of $\mathbb{T}$}, denoted by $\text{Aut}(\mathbb{T})$. Every automorphism of $\mathbb{T}$ is continuous and thus a homeomorphism of the torus, cf. \cite[Theorem 4.4.1]{gunter2001}. With respect to the compact-open topology, the group $\text{Aut}(\mathbb{T})$ is a Lie group of dimension at most 6, cf.  \cite[Theorem 1.1]{brendan2017a}. 

The automorphism group $\text{\normalfont Aut}(\mathbb{T})$ has two distinguished normal subgroups, the \textit{kernels} $T^\pm$  consisting of all automorphisms of $\mathbb{T}$ that fix every $(\pm)$-parallel class. 
For convenience, we  refer to these two subgroups as the \textit{kernels $T^\pm$  of the plane $\mathbb{T}$}.   
	The connected component of $T^\pm$ of a flat Minkowski plane is isomorphic to $\text{\normalfont PSL}(2,\mathbb{R})$, $\text{\normalfont L}_2,$ $\text{\normalfont SO}(2,\mathbb{R})$, $\mathbb{R}$,   or the trivial group $\{id\}$, cf.  \cite[Proposition 4.4.9]{gunter2001}. The same result holds for toroidal circle planes, since their automorphism groups are Lie groups. 
	

Another fact that we frequently use is that automorphisms fixing three pairwise non-parallel points have order at most 2. In particular, if such an automorphism takes $(+)$-parallel classes to $(+)$-parallel classes, then it is the identity map. A proof can be found in \cite[Lemma 4.4.2]{gunter2001}.

 \subsection{Some results on transformation groups}
 
The following theorem describes all possible transitive and effective actions of transformation
groups on 1-manifolds. We refer to this result as Brouwer's Theorem throughout
this paper.
In   some sources, a weaker version (for Lie groups) is stated, without a name, as a consequence of a result by Sophus Lie on Lie algebras, cf. \cite{navas2011},  \cite[p.348]{ghys2001}  and \cite[Theorem 2.1, p.218]{onishchik1993}. According to \cite[96.30]{salzmann1995}, this theorem is proved by Brouwer \cite{brouwer1909}, and a sketch of proof is provided there.
 
\begin{theorem}[Brouwer's Theorem]  \label{brouwer}
	Let $G$ be a locally compact, connected, effective and transitive transformation group on a connected 1-dimensional manifold M. Then G has dimension at most 3. 
	
	(a) If $M \cong \mathbb{S}^1$, then $G$ is isomorphic and acts equivalently to the rotation group $\text{\normalfont SO}(2,\mathbb{R})$ or a finite covering group $\text{\normalfont PSL}^{(k)}(2,\mathbb{R})$ of the projective group $\text{\normalfont PSL}(2,\mathbb{R})$. 
	
	(b) If $M \cong \mathbb{R}$, then $G$ is isomorphic and acts equivalently to $\mathbb{R}$, the connected component $\text{\normalfont L}_2$ of the affine group of $\mathbb{R}$, or the simply connected covering group $\widetilde{\text{\normalfont PSL}(2,\mathbb{R})}$ of  $\text{\normalfont PSL}(2,\mathbb{R})$. 
\end{theorem}

Often,  groups of automorphisms are assumed to be connected, so that they map $(\pm)$-parallel classes to $(\pm)$-parallel classes.  Consequently, they have induced actions on $\mathcal{G}^\pm$, which are homeomorphic to $\mathbb{S}_1$. It is then important to know how these
groups act on $\mathbb{S}_1$ and $\mathbb{R}$  (the only connected 1-manifolds). The following is helpful for this purpose.

\begin{lemma}[cf. \cite{salzmann1995} 96.29]\label{salzmann9629} 
	Consider a connected group $\Gamma$ acting effectively on $\mathbb{R}$ or $\mathbb{S}^1$.
	
	(a) If $\Gamma$ has no fixed point, then $\Gamma$ is transitive.
	
	(b) Any non-trivial compact subgroup of $\Gamma$ acts freely on $\mathbb{S}^1$; it cannot act on $\mathbb{R}$. 
	
	(c) If $\Gamma$ is compact, then $\Gamma =\{id\}$, or $\Gamma \cong \text{\normalfont SO}(2,\mathbb{R})$ and $\Gamma$ is sharply transitive on $\mathbb{S}^1$. 
\end{lemma}
%
\begin{lemma}[cf. \cite{salzmann1995} 93.12] \label{salzmann9312}  If $\Delta$ is a closed subgroup of the locally compact, connected group $\Gamma$, and if $\dim \Delta = \dim \Gamma < \infty$, then $\Delta =\Gamma$. 
\end{lemma} 
Proofs of the main theorems rely on arguments with the dimension of orbits and stabilisers, which
are based on the following.
\begin{lemma}[cf. \cite{salzmann1995} 96.10, \cite{gunter2001} Theorem A.2.3.6] \label{dimensionformula} If the Lie group $G$ acts on a manifold $M$, then
	$$
	\dim G= \dim G_p +\dim G(p),
	$$
	where $G_p$ and $G(p)$ are the stabiliser and orbit, respectively, of the point $p \in M$. 
\end{lemma}
We will usually refer to Lemma \ref{dimensionformula} as the dimension formula.

\section{Proof of Theorem \ref{almostsimplePSL}}

Let $\mathbb{T}$ be a toroidal circle plane with automorphism group 
$\text{\normalfont Aut}(\mathbb{T})$. 
We denote the connected component of $\text{\normalfont Aut}(\mathbb{T})$ by $\Gamma$ and let $S$  be a non-trivial almost simple connected Lie subgroup of $\Gamma$. Let $K^\pm$ be the kernel of the action of $S$ on $\mathcal{G}^\pm$. 

The overall structure of the proof of Theorem \ref{almostsimplePSL} is  as follows. We first show in Lemma \ref{PSLlocally} that $S$ is locally isomorphic to $\text{\normalfont PSL}(2,\mathbb{R})$. 
This implies $S$ is isomorphic to either the universal covering group $\widetilde{\text{\normalfont PSL}(2,\mathbb{R})}$ or a finite covering group $\text{\normalfont PSL}^{(k)}(2,\mathbb{R})$ of $\text{\normalfont PSL}(2,\mathbb{R})$. 
In particular, the centre $Z(S)$ is cyclic. 
Next, we prove that $S$ is in fact isomorphic to $\text{\normalfont PSL}(2,\mathbb{R})$  by showing that the centre $Z(S)$ is trivial in Lemma \ref{trivialcentre}. Finally, in Lemma \ref{actionalmostsimple}, we apply Brouwer's Theorem (Theorem \ref{brouwer}) to determine possible actions of $S$ on the torus. 

We start with  the following observation.

\begin{lemma} \label{SO2Rsubgroup} If $S$ contains a subgroup $H$ isomorphic to $\text{\normalfont SO}(2,\mathbb{R})$, then $S$ acts transitively on at least one of $\mathcal{G}^\pm$.
	\begin{proof} Since $T^+ \cap T^- = \{id \}$, the subgroup $H$  cannot be contained in both $T^\pm$. Without loss of generality, we assume $H \not \subset T^+$. This means  $H$ cannot act trivially, and therefore, acts transitively on $\mathcal{G}^+ \cong \mathbb{S}^1$ (cf. Lemma \ref{salzmann9629}). 
	\end{proof}
\end{lemma}

\begin{lemma} \label{PSLlocally} $S$ is locally isomorphic to $\text{\normalfont PSL}(2,\mathbb{R})$.
	\begin{proof}
		We first note that $S$ cannot have dimension 6, because this implies $\mathbb{T}$ is the classical flat Minkowski plane (cf. \cite[Theorem 2]{brendan2017a}) and $S \cong \text{\normalfont PSL}(2,\mathbb{R}) \times \text{\normalfont PSL}(2,\mathbb{R})$, which is not almost simple.
		From the list of almost simple groups of low dimensions (cf. \cite[Theorem A2.2.6]{gunter2001} or \cite[\text{94.33}]{salzmann1995}), $S$ is locally isomorphic to either $\text{\normalfont SO}(3,\mathbb{R})$ or $\text{\normalfont PSL}(2,\mathbb{R})$.

		Suppose for a contradiction that $S$ is locally isomorphic to $\text{\normalfont SO}(3,\mathbb{R})$, that is, $S$ is isomorphic to either  $\text{\normalfont SO}(3,\mathbb{R})$ or  $\text{\normalfont SU}(2,\mathbb{C})$. Then $S$ is 3-dimensional, compact and contains a subgroup isomorphic to $\text{\normalfont SO}(2,\mathbb{R})$. By Lemma \ref{SO2Rsubgroup}, we may assume $S$ is transitive on $\mathcal{G}^+$. By Brouwer's Theorem and the compactness of $S$, the factor group $S/K^+$ is isomorphic to  $\text{\normalfont SO}(2,\mathbb{R})$. But this implies $\dim K^+ = 2$, which contradicts the assumption $S$ is almost simple. 
	\end{proof}
\end{lemma}


Let $\kappa$ be a generator of $Z(S)$. 

\begin{lemma} \label{PSLpart1} $Z(S)$ is contained in at least one of $K^\pm$.
	\begin{proof}  
		In parts 1) to 5) we show that  $\kappa(p) \parallel p$ for every point $p \in \mathcal{P}$. The final step, part 6), yields the lemma. Suppose for a contradiction that there exists a point $p$ such that $\kappa(p)$ is not parallel to  $p$.

		\begin{enumerate}[label=\arabic*),leftmargin=0pt,itemindent=*]
		
			\item We show  $\dim S_{[p]_\pm} =2$. From Lemma \ref{dimensionformula},
			$$
			3= \dim S= \dim S_{[p]_+} + \dim S([p]_+).
			$$
			Since $\dim S([p]_+)$ is either $0$ or $1$, $\dim S_{[p]_+}$ is either $3$ or $2$.  If $\dim S_{[p]_+} =3$, then by Lemma \ref{salzmann9312}, $S= S_{[p]_+}$. This cannot be true however, since $\kappa$ does not fix $[p]_+$.  Hence $\dim S_{[p]_+} =2$.
			Similarly, $\dim S_{[p]_-} =2$.
			
			\item  We show that 	 $S_{[p]_\pm}$  fixes at least one point on $[p]_\pm$.  If 	$S_{[p]_+}$ is transitive on $[p]_+$, then by Brouwer's Theorem $S_{[p]_+}/ K^- \cong \text{\normalfont SO}(2,\mathbb{R})$. Since $\dim S_{[p]_+} =2$, $\dim K^- =1$. But this contradicts the fact that $S$ is almost simple. Hence $S_{[p]_+}$  is not transitive on $[p]_+$. By Lemma \ref{salzmann9629}, $S_{[p]_+}$ fixes at least one point on $[p]_+$. The same argument holds for $S_{[p]_-}$ on $[p]_-$.

			\item We show that either $S_p=S_{[p]_-}$ or $S_p=S_{[p]_+}$. Following part 2), let $q$ be a point such that $[q]_\mp$ is fixed by $S_{[p]_\pm}$. Since  $S_p \le  S_{[p]_+} \cap S_{[p]_-}$, $S_p=S_{p,q}$. Because $\kappa \in Z(S)$, we also have $S_p=S_{p,\kappa(p)}$, and so $S_p=S_{p,q,\kappa(p)}$. If $q$ is nonparallel to $p$ and $\kappa(p)$, then $S_p$ is trivial and $\dim S_p = 0$. From the dimension formula, we get $\dim S(p)= \dim S - \dim S_p = 3$, which cannot be true. Hence either $q \parallel p$ or $q \parallel \kappa(p)$. 
			
			If $q \parallel_+ p$, then $[q]_+=[p]_+$, so that $S_{[p]_-}$ fixes $[p]_+$. If  $q \parallel_+ \kappa(p)$, then $S_{[p]_-}$ fixes $\kappa([p]_+)$, and since $\kappa$ commutes with every element of $S$, we see that $S_{[p]_-}$ also fixes $[p]_+$.  In both cases, $S_p=S_{[p]_-}$.
			
			In the cases $q \parallel_- p$ and $q \parallel_- \kappa(p)$, we obtain $S_p=S_{[p]_+}$ in a similar manner. 
			
			\item We show that there exists $r \in S(p)$ such that $r \not \parallel p$ and $r \not \parallel \kappa(p)$. From part 1) and 3),  $\dim S_p=2$, so that $\dim S(p)=1$.  Also from part 3), we can assume $S_p=S_{[p]_-}$, and so $S_{\kappa(p)} =S_{\kappa([p]_-)}$. 
			
			Suppose for a contradiction that   for  every $\alpha \in S$, either  $\alpha(p) \parallel p$ or $\alpha(p) \parallel \kappa(p)$.   
			From the assumption $S_p=S_{[p]_-}$, we see that $\alpha(p)=p$ if $\alpha(p) \parallel_- p$.   
			Similarly, if $\alpha(p) \parallel_- \kappa(p)$, then  $\alpha$ fixes $\kappa(p)$. 
			It follows that  the orbit $S(p)$ consists of points in $[p]_+ \cup [\kappa(p)]_+$. But this is impossible, since $S(p)$ is connected. 
				\item Let $r$ be as in part 4). From the dimension formula, we have 
			\begin{align*}
			\dim S_p 	&=  \dim S_{p,r}+\dim S_p(r)  \le \dim S_{p,r, \kappa(p)}+1=1, 
			\end{align*} 
			a contradiction. This proves that  $\kappa(p) \parallel p$ for every point $p \in \mathcal{P}$.
			\item	 Suppose that there are two nonparallel points $p,q$ such that $\kappa(p) \ne p, \kappa(q) \ne q$, $\kappa(p) \parallel_+p$, $\kappa(q) \parallel_- q$. Let $r= [p]_- \cap [q]_+$. Then $\kappa(r) \parallel_- \kappa(p)$ and $\kappa(r) \parallel_+ \kappa(q)$, compare Figure \ref{fig:centreinkernel}.
			Since $\kappa(p) \ne p$ and  $\kappa(q) \ne q$, we see that $\kappa(r)$ is nonparallel to $r$, a contradiction. 
			This shows that  $Z(S)$ is contained in at least one of $K^\pm$.  \qedhere

			\begin{figure}[h]
				\begin{center}
					\begin{tikzpicture}[scale=0.5]
					
					\draw [fill] (-4,-2) circle [radius=0.1];
					\node[below left] at (-4,-2) {$p$};
					\coordinate (d_1) at (-4,-2);
					
					\draw [fill] (0,0) circle [radius=0.1];
					\node[below left] at (0,0) {$\kappa(r)$};
					\coordinate (d_2) at (0,0);
					
					\draw [fill] (0,7) circle [radius=0.1];
					\node[above left] at (0,7) {$\kappa(q)$};
					\coordinate (d_4) at (0,7);

					\draw [fill] (6,-2) circle [radius=0.1];
					\node[below right] at (6,-2) {$r$};
					\coordinate (d_5) at (6,0);

					\draw [fill] (6,7) circle [radius=0.1];
					\node[above right] at (6,7) {$q$};
					\coordinate (d_3) at (6,7);
					
					\draw [fill] (-4,0) circle [radius=0.1];
					\node[below left] at (-4,0) {$\kappa(p)$};
					\coordinate (kp) at (-4,0);
					
					\coordinate (z_1) at (-4,7);
					\coordinate (z_2) at (-4,0);
					\coordinate (z_3) at (0,-2);
					\coordinate (z_4) at (6,-2);5
					
					%
					\draw   (-4,8) -- (-4,-3);
					\draw    (-0,8) -- (-0,-3);
					\draw   (6,8) -- (6,-3);
					\draw    (-5,7) -- (7,7);
					\draw   (-5,0) -- (7,0);
					\draw   (-5,-2) -- (7,-2);
					
					\end{tikzpicture}
				\end{center}
				\caption{} \label{fig:centreinkernel}
			\end{figure}
			 
		\end{enumerate} 
		
	\end{proof}	
\end{lemma}

\begin{lemma} \label{trivialcentre} $Z(S)$ is trivial.
In particular, $S \cong \text{\normalfont PSL}(2,\mathbb{R})$.
	
	\begin{proof}    
		
		Since $K^\pm$ is a normal subgroup of $S$, the dimension of $K^\pm$ is either $0$ or $3$.  
		If either $\dim K^\pm =3$, then by Lemma \ref{salzmann9312}, the list of possible connected groups in a kernel (cf. \cite[Proposition 4.4.9]{gunter2001}), and the fact that $S$ is almost simple, it follows that $S =K^\pm \cong \text{\normalfont PSL}(2,\mathbb{R})$ and thus $Z(S)$ is  trivial.

		In the remainder  of the proof, we deal with the case that both	$\dim K^\pm = 0$. Following Lemma \ref{PSLpart1}, without loss of generality, we assume $Z(S) \le K^+$.   As a discrete normal subgroup of $S$, $K^+ \le Z(S)$ and hence $Z(S) = K^+$. Also,  $K^- \le Z(S)$ and since $K^-\cap K^+=\{id\}$, we see that $K^-$ is trivial. 
		
		Suppose for a contradiction that $Z(S)$ is non-trivial. 
		

		\begin{enumerate}[label=\arabic*),leftmargin=0pt,itemindent=*] 
			
			\item	We show that the action of $S/K^+$ on $\mathcal{G}^+$ is equivalent to  the standard action of $\text{\normalfont PSL}(2,\mathbb{R})$ on $\mathbb{S}^1$. We first note that $S/K^+ =S/Z(S) \cong \text{\normalfont PSL}(2,\mathbb{R})$ so that $S/K^+$ contains a subgroup $H$ isomorphic to $\text{\normalfont SO}(2,\mathbb{R})$. By Lemma \ref{salzmann9629}, $H$   acts  transitively on $\mathcal{G}^+$, as the preimage of $H$ in $S$ is not   a subgroup of $K^+$. Hence $S/K^+$ is transitive on  $\mathcal{G}^+$ and its action is derived from Brouwer's Theorem. 
		\item   Since $S \ne K^-$, there exists a point $x$ such that $\dim S([x]_-)=1$.  From the dimension formula,  $\dim S_{[x]_-} =2$. 
Similar to part 2) in the proof of Lemma \ref{PSLpart1},  $S_{[x]_-}$ fixes a point $p$ on $[x]_-$. This implies $S_{[x]_-} \le S_p$, and in particular,
$
2 \le \dim S_p \le 3. 
$ 
Furthermore, $\dim S_p \ne 3$, otherwise from Lemma \ref{salzmann9312}, $S_p=S$, a contradiction to the transitivity of $S/K^+$ on $\mathcal{G}^+$.  Hence $\dim S_p =2$ and $\dim S(p)=1$.  For the rest of the proof, we will fix such a  point $p$.

\item 
Denote the connected component of $S_p$ by $S_p^1$.    

If $S_p^1$ is not transitive on $\mathcal{G}^+\backslash\{[p]_+\}$, then it fixes a point $q \in [\kappa(p)]_-$ nonparallel to $p$. 
From the dimension formula, for a point $r \in S(p)$ nonparallel to $p$ and $q$, we have
$$
2= \dim S^1_p = \dim  S^1_{p,q}= \dim S^1_{p,q,r}+ \dim S^1_{p,q}(r),
$$
which implies $ \dim S^1_{p,q}(r) = 2$, contradicting the fact that $\dim S(p)=1$.
Hence $S_p^1$ is transitive on $\mathcal{G}^+\backslash\{[p]_+\}$. 
Since $S_p^1 \cap K^+$ is trivial, Brouwer's Theorem implies that  $S_p^1$ is isomorphic and acts     equivalently to $\text{\normalfont L}_2$ on $\mathcal{G}^+\backslash\{[p]_+\}$.

Let $R \cong \mathbb{R}$ be a 1-dimensional orbit of $S_p^1$  on $[p]_+$. 
Let $R^-$ be the kernel of the action of $S_p^1$ on $R$. 
By Brouwer's Theorem,  either $S^1_{p}/R^- \cong \mathbb{R}$ or $S^1_{p}/R^- \cong \text{\normalfont L}_2$.  
If  $S^1_{p}/R^- \cong \mathbb{R}$, then   $R^-$ is isomorphic to $\mathbb{R}$, and, as a 1-dimensional normal subgroup of $S^1_{p} \cong \text{\normalfont L}_2$, acts transitively on $\mathcal{G}^+\backslash\{[p]_+\}$. 
However, this implies $[x]_- \backslash \{x\} \in S(p)$  for each    $x\in R$, which cannot occur, since $\dim S(p)=1$. 
Hence, $S^1_{p}/R^- \cong \text{\normalfont L}_2$. In particular,  $R^-$ is trivial and   $S_p^1$    acts equivalently to $\text{\normalfont L}_2$ on $R$.

It follows that $S_p^1$ has a 1-dimensional orbit $Q \cong \mathbb{R}$ on $R \times \mathcal{G}^+\backslash\{[p]_+\}$, which corresponds to  the diagonal under the identification $R \times \mathcal{G}^+\backslash\{[p]_+\} \cong \mathbb{R}^2$, cf. Figure \ref{fig:trivialcentre1}. 
Note that any two distinct points in $Q$ are non-parallel so that $S_p^1$ acts equivalently to $\text{\normalfont L}_2$ on $Q$.


			\begin{figure}[h] 
							\centering
					\begin{tikzpicture}[scale=0.65]

					\draw [fill] (0,0) circle [radius=0.1];
					\coordinate (d_2) at (0,0);
					
					\draw [fill] (0,7) circle [radius=0.1];
					\node[above left] at (0,7) {$p$};
					\coordinate (d_4) at (0,7);

					\draw [fill] (7,0) circle [radius=0.1];
					\coordinate (d_5) at (7,0);

					\draw [fill] (7,7) circle [radius=0.1];
					\node[above right] at (7,7) {$p$};
					\coordinate (d_3) at (7,7);
%
%
					\draw [fill] (0,6) circle [radius=0.1];
					\coordinate (kp1) at (0,6);

					\draw [fill] (7,6) circle [radius=0.1];
					\coordinate (kp2) at (7,6);
					
					\draw  (d_4) -- (d_3);
					\draw  (d_5) -- (kp2);
					\draw  [dashed] (kp2) -- (d_3);
					\draw  [dashed] (d_4) -- (kp1);
					
					\draw  (kp1) -- (d_2);
					\draw  (d_2) -- (d_5);

					\draw  (d_2) -- (kp2); 
					\draw (kp1) -- (kp2); 
					
					\node[right] at (3.8,3) {$Q$};  
					\node[left] at (0,3) {$R$};
					
						\draw [fill] (0,-1) circle [radius=0.1];
					\node[below left] at (0,-1) {$p$};  
					\draw [fill] (7,-1) circle [radius=0.1];
					\node[below right] at (7,-1) {$p$}; 
					\draw  [dashed] (0,-1) -- (0,0);
					\draw  [dashed] (7,-1) -- (7,0);
					\draw  (0,-1) -- (7,-1);
					\end{tikzpicture}  
				\caption{} \label{fig:trivialcentre1}
			\end{figure}

\item
 
For a given point $x \in Q$, let  $\mathcal{C}_{p,x}$ be the set of circles going through $p$ and $x$. Let  $\mathcal{C}_{p,x}^+$ and $\mathcal{C}_{p,x}^-$ be the subsets of $\mathcal{C}_{p,x}$ consisting of circles described by orientation-preserving and orientation-reversing homeomorphisms of $\mathbb{S}^1$ to itself. 
Define $$ \phi:  
 [\kappa(p)]_- \backslash \{\kappa(p),  [\kappa(p)]_- \cap [x]_+\} \rightarrow \mathcal{C}_{p,x}  :
 y \mapsto \alpha(p,x,y).$$

The map $\phi$ is induced from the operation Joining, which is a homeomorphism, cf. \cite{duythesis}.   
In particular, $\phi$ is also a homeomorphism and   maps each connected component of its domain onto  one of $\mathcal{C}_{p,x}^\pm$. 

\item  We now make two observations 
on $S_{p,x}^1$,   the connected component of the 2-point stabilizer $S_{p,x}$, for each $x \in Q$. 
Firstly, from part 3),  $S_{p,x}^1$ is isomorphic to $\mathbb{R}$ and acts sharply transitively on each  connected component of $Q \backslash \{x\}$.  

Secondly,  $S_{p,x}^1$ is also sharply transitive on each of $\mathcal{C}_{p,x}^\pm$.  
This comes from the action of $S_{p}^1$ on  the fixed parallel class $[\kappa(p)]_-$ and the identification of  points on $[\kappa(p)]_-$ with circles in $\mathcal{C}_{p,x}^\pm$ described in part 4).
 
 	\item   Let $q \in Q$ be a point  and $U,V$ be the two connected components of $Q \backslash \{q\}$.   
 	Let $u \in U, v \in V$ and let $A,B$ be the circles $A:=\alpha(p,q,u)$ and $B:=\alpha(p,q,v)$, cf. Figure \ref{fig:trivialcentre2}.

By part 3),   $S_p^1$ is 2-set transitive (2-homogeneous) on $Q$, and so there exists $\gamma \in S_p^1$ such that $\gamma(\{v,q\})=\{q,u\}$.  In particular, $\gamma(A) =B$.  
Since $S_p^1$ is connected,  both $A,B$ belong to either $\mathcal{C}_{p,q}^+$ or $\mathcal{C}_{p,q}^-$. Assume $A,B \in \mathcal{C}_{p,q}^+$.
			
			\begin{figure}[h]
				\begin{center}
					\begin{tikzpicture}[scale=0.65]
					\draw [fill] (0,0) circle [radius=0.1];
					\coordinate (d_2) at (0,0);
					
					\draw [fill] (0,7) circle [radius=0.1];
					\node[above left] at (0,7) {$p$};
					\coordinate (d_4) at (0,7);

					\draw [fill] (7,0) circle [radius=0.1];
					\coordinate (d_5) at (7,0);

					\draw [fill] (7,7) circle [radius=0.1];
					\node[above right] at (7,7) {$p$};
					\coordinate (d_3) at (7,7);
					
					\draw [fill] (2.5,2.13) circle [radius=0.1];
					\node[above left] at (2.5,2.13) {$q$};
					\coordinate (q) at (2.5,2.13);
					
					\draw [fill] (0,6) circle [radius=0.1];
					\coordinate (kp1) at (0,6);

					\draw [fill] (7,6) circle [radius=0.1];
					\coordinate (kp2) at (7,6);
					
					\coordinate (r2) at (6,6);
					
					\coordinate (s2) at (5.5,6);
					
					\draw  (d_4) -- (d_3);
					\draw  (d_5) -- (kp2);
					\draw  [dashed]  (kp2) -- (d_3);
					\draw  [dashed]  (d_4) -- (kp1);
					\draw  (kp1) -- (d_2);
					\draw  (d_2) -- (d_5);3

					\draw  (d_2) -- (kp2); 
					\draw (kp1) -- (kp2);

					\draw [densely dashed] (0,-1)  to [out=75,in=180] (q) to [out=0,in=220] (d_3)  ; 
					
						\coordinate (v) at (1.22,1.05);
						
						\draw [fill] (v) circle [radius=0.1];
						\node[below right] at (v) {$v$};
					
					\draw [thick] (0,-1)  to [out=45,in=225] (v) to [out=45,in=270] (q) to [out=90,in=200] (d_3)  ;

					\draw [fill] (5.27,4.49) circle [radius=0.1];
					\node[below right] at (5.27,4.49) {$u$};
					
					\node[below right] at (4.2,2.9) {$A$};
					\node[left] at (3.1,4.6) {$B$}; 
					
					\draw [fill] (0,-1) circle [radius=0.1];
					\node[below left] at (0,-1) {$p$};  
					\draw [fill] (7,-1) circle [radius=0.1];
					\node[below right] at (7,-1) {$p$}; 
					\draw  [dashed] (0,-1) -- (0,0);
					\draw  [dashed] (7,-1) -- (7,0);
					\draw  (0,-1) -- (7,-1);
					
						\node[left] at (0,3) {$R$}; 
					\end{tikzpicture}
				\end{center}
				\caption{} \label{fig:trivialcentre2}
			\end{figure}

			From part 5), there exists a unique $\sigma \in S_{p,q}^1$ such that $\sigma(B)=A$.  This implies $\sigma(v) \in A$, and from the action of $S_{p,q}^1$ on $Q$ in part 5), it follows that $ \sigma(v)  \in V$.  Since $\sigma(v)$ and $q$ are in the same connected component of $Q \backslash \{ u\}$, there exists a unique $\tau \in S_{p,u}^1$ such that $\tau(q)=\sigma(v)$. In other words,  $\tau$ fixes $p,u \in A$ and maps $q \in A$ to $\sigma(v) \in A$, so $\tau$ must fix $A$.  
			This implies
			$\tau = id$. 
			But this is impossible, because $\sigma(v)  \ne q$. This completes the proof. \qedhere

		\end{enumerate}
		
	\end{proof}
\end{lemma} 

We now determine the action of $S$ on the torus.

\begin{lemma} \label{actionalmostsimple} If $S \cong \text{\normalfont PSL}(2,\mathbb{R})$,   then exactly one of the following occurs.
	\begin{enumerate}
		\item Either $S$ fixes every $(+)$-parallel class or $S$ fixes every $(-)$-parallel class. In both cases, $\mathbb{T}$ is isomorphic to a half-classical Minkowski plane $\mathcal{M}(f,id)$, where $f$ is an orientation-preserving homeomorphism of $\mathbb{S}^1$. 
		
		\item  $S$ acts diagonally on the point set. The diagonal $D$ is a circle of $\mathbb{T}$, and $S$ fixes $D$. 
	\end{enumerate}
	\begin{proof}  Since $K^\pm$ are normal subgroups of the simple group $S$,  either $K^\pm = S$ or $K^\pm= \{ id \}$. 
		If  $K^+ = S$, then $S$ fixes every $(+)$-parallel class, and likewise when $K^- = S$. In both cases, one of the kernels of $\mathbb{T}$ is 3-dimensional. From \cite[Theorem 1.2]{brendan2017a} and \cite[Theorem 4.4.10]{gunter2001}, $\mathbb{T}$ is determined.
		
		We now consider the case that both $  K^\pm=\{ id \}$. Since $S$ contains   a subgroup $H \cong \text{\normalfont SO}(2,\mathbb{R})$ which is not contained in $K^\pm$, by Lemma \ref{salzmann9629}, $H$ and thus $S$ acts transitively on $\mathcal{G}^\pm$. 
		Hence  $S$ acts equivalently to $\text{\normalfont PSL}(2,\mathbb{R})$ on both $\mathcal{G}^\pm$. It follows that $S$ has two orbits: the diagonal $D \cong \mathbb{S}^1$, in a suitable coordinate system,  and its complement $\mathcal{P} \backslash D$.

		It only remains to show that $D$ is a circle of $\mathbb{T}$.

		\begin{enumerate}[label=\arabic*),leftmargin=0pt,itemindent=*]
			
			\item As seen in the proof of Lemma \ref{trivialcentre}, for any pair of distinct points $p,q \in D$, the stabiliser $S_{p,q}$ is isomorphic to $\mathbb{R}$.  The orbits of points in $\mathcal{P} \backslash \{[p],[q] \}$ under $S_{p,q}$ are 1-dimensional and partition $\mathcal{P} \backslash \{[p],[q] \}$. Under a suitable coordinate system with $p=(\infty,\infty)$ and $q= (0,0)$,  these orbits can be represented as sets of the form $\{ (x,ax)|x>0 \}$ or $\{ (x,ax)|x<0 \}$, 
			where each $a\ne 0$ determines such an orbit. 
			
			\item Let  $p,q,r \in D$ and $C \coloneqq \alpha(p,q,r)$.
			We claim that   there exists an  orbit   $O \coloneqq S_{x,y}(\xi)$,  where $x,y \in \{p,q,r\},  x\ne y$, for some $\xi \in \mathcal{P}$, intersecting $C$ at least two times.   
			If $C$ intersects $D$ at an additional point $v$, then depending on the position of $v$, the orbit $O$ can be chosen as  $S_{p,q}(r)$, $S_{p,r}(q)$, or $S_{q,r}(p)$, cf. Figure \ref{fig:pslaction1}.  
			 
			Otherwise, at least one intersection of $C$ and $D$ is not transversal.   
			By changing the roles of $p,q,r$ if necessary, we can assume this intersection is $r$, cf. Figure \ref{fig:pslaction2}. 
			 Coordinatise the plane as in part 1) so that the orbits of points under $S_{p,q}$ are rays emanating from $q$. Let $C_0$ be the connected component of $C\backslash \{p,q\}$ that contains $r$.   Let $C_1,C_2$ be the two connected components of $C_0 \backslash \{r\}$.  
			Since the orbits $S_{p,q}(\xi)$ depend continuously on $\xi$ and approach $S_{p,q}(r) \subseteq D$ as $\xi$ tends to $r$,
			the intermediate value theorem implies that there exists an orbit that meets $C$ in at least two points, one in $C_1$ and one in $C_2$.  
		 
			\begin{figure}[h] 
				\centering
				\begin{subfigure}[h]{0.45\textwidth}
					\centering
					\begin{tikzpicture}[scale=0.78]
					\draw [fill] (0,0) circle [radius=0.1];
					\node[below left] at (0,0) {$p$};
					\coordinate (d_2) at (0,0);
					
					\draw [fill] (0,7) circle [radius=0.1];
					\node[above left] at (0,7) {$p$};
					\coordinate (d_4) at (0,7);

					\draw [fill] (7,0) circle [radius=0.1];
					\node[below right] at (7,0) {$p$};
					\coordinate (d_5) at (7,0);

					\draw [fill] (7,7) circle [radius=0.1];
					\node[above right] at (7,7) {$p$};
					\coordinate (d_3) at (7,7);
					
					\draw [fill] (2,2) circle [radius=0.1];
					\node[below right] at (2,2) {$q$};
					\coordinate (q) at (2,2);
					
					\draw [fill] (4,4) circle [radius=0.1];
					\node[left] at (4,4) {$r$};
					\coordinate (r) at (4,4);
					
					\draw [fill] (6,6) circle [radius=0.1];
					\node[above left] at (6,6) {$v$};
					\coordinate (v) at (6,6);
					
					\draw  (d_4) -- (d_3);
					\draw  (d_5) -- (d_3);
					\draw  (d_4) -- (d_2);
					\draw  (d_2) -- (d_5);

					\draw  (d_2) -- (q); 
					\draw  [dashed] (d_3) -- (q); 
					
					\draw [thick] (d_2)  to [out=70,in=180] (q) to [out=0,in=270] (r) to  [out=90,in=180] (v) to  [out=0,in=270] (d_3) ; 
					%
					%
					%
					\node[right] at (5,5) {$O$};
					\node[left] at (4.7,5.6) {$C$};
					
					\end{tikzpicture}
					\caption{} \label{fig:pslaction1}
				\end{subfigure}
				\hfill
				\begin{subfigure}[h]{0.45\textwidth}
					\centering
					\begin{tikzpicture}[scale=0.78]
					\draw [fill] (0,0) circle [radius=0.1];
					\node[below left] at (0,0) {$p$};
					\coordinate (d_2) at (0,0);
					
					\draw [fill] (0,7) circle [radius=0.1];
					\node[above left] at (0,7) {$p$};
					\coordinate (d_4) at (0,7);

					\draw [fill] (7,0) circle [radius=0.1];
					\node[below right] at (7,0) {$p$};
					\coordinate (d_5) at (7,0);

					\draw [fill] (7,7) circle [radius=0.1];
					\node[above right] at (7,7) {$p$};
					\coordinate (d_3) at (7,7);
					
					\draw [fill] (2,2) circle [radius=0.1];
					\node[below right] at (2,2) {$q$};
					\coordinate (q) at (2,2);
					
					\draw [fill] (4,4) circle [radius=0.1];
					\node[left] at (4,4) {$r$};
					\coordinate (r) at (4,4);
					
					\draw [fill] (4.47,3.5) circle [radius=0.1];
					\node[below right] at (4.4,3.5) {$\xi$}; 
					
					\draw [fill] (3.,2.5) circle [radius=0.1];
					\node[below right] at (3,2.5) {$s$};

					\draw [fill] (6.1,5.4) circle [radius=0.1];
					\node[below right] at (6.1,5.4) {$t$}; 
					
					\draw  (d_4) -- (d_3);
					\draw  (d_5) -- (d_3);
					\draw  (d_4) -- (d_2);
					\draw  (d_2) -- (d_5);

					\draw  (d_2) -- (d_3); 
					
					\draw [thick] (d_2)  to [out=70,in=180] (q) to [out=0,in=225] (r)   to  [out=45,in=270] (d_3) ; 
					
					\draw [dashed]  (q) to [out=25,in=245]  (d_3) ;
					
					\node[right] at (5.2,4) {$O$};
					\node[left] at (1.3,2) {$C$};
					
					\end{tikzpicture}
					\caption{} \label{fig:pslaction2}
				\end{subfigure}
				\caption{} \label{fig:pslaction}
			\end{figure}

			\item  
Following part 2), without loss of generality, assume $O \coloneqq S_{p,q}(\xi)$ for some $\xi \in \mathcal{P}$. 
			 Let $R$ be the connected component of $D \backslash \{p,q\}$ that contains $r$. We show that $R$ intersects $C$  infinitely many times.
			 
			Let $s,t \in O \cap C,$ $s \ne t$. Then there exists $\sigma \in S_{p,q}$ such that  $\sigma(s)=t$. This implies $\sigma(C)=C$. But since $s \ne t$, we have  $\sigma \ne id$, so that $\sigma(r) \ne r$ and $\sigma(r) \in R \cap C$. 
			Since $\sigma$ has infinite order,   the points $\sigma^n(r)$ are distinct and belong to $R \cap C$. 
			
			\item By repeating the arguments in part 2) and 3) for different pairs of points in $D \backslash \{ p\}$,  it follows that the set $D \cap C$ is dense on $D$. From the compactness of $D$ and $C$, we have $D \subseteq C$. As there is no proper subset of $\mathbb{S}^1$ homeomorphic to $\mathbb{S}^1$, it follows that $D=C$. This completes the proof. \qedhere
		\end{enumerate}
	\end{proof}
\end{lemma}

 This proves Theorem \ref{almostsimplePSL}. We also obtain the following results. 
\begin{corollary}\label{solvable}
	Let $\Sigma$ be a 3-dimensional connected group of automorphisms of a toroidal circle plane. Then  $\Sigma$ is either solvable or isomorphic to $\text{\normalfont PSL}(2,\mathbb{R})$. 
%
\end{corollary}

\begin{corollary}\label{semisimple}
	Let $\Sigma$ be a connected semi-simple group of automorphisms of a toroidal circle plane. Then  $\Sigma$ is isomorphic to either  $\text{\normalfont PSL}(2,\mathbb{R}) \times \text{\normalfont PSL}(2,\mathbb{R})$ or $\text{\normalfont PSL}(2,\mathbb{R})$. 
\end{corollary}

\section{Proof of Theorem \ref{3dfull}} 

In this section, let $\Sigma$ be a 3-dimensional connected group of automorphisms of a toroidal circle plane $\mathbb{T}$ with kernels $\Delta^\pm$  on $\mathcal{G}^\pm$. 
We prove Theorem \ref{3dfull}  in three parts via Lemmas \ref{3dpc}, \ref{3dfixedelements} and \ref{3dfixepoint}.  
\begin{lemma} \label{3dpc} Exactly one of the following occurs.
	
	\begin{enumerate}
		\item 		$\Sigma$ fixes at least one parallel class.
		
		\item  $\Sigma$ fixes no parallel classes but fixes and acts transitively on exactly one circle. In this case $\Sigma \cong \text{\normalfont \text{\normalfont PSL}}(2,\mathbb{R})$ and, under suitable coordinates, it acts diagonally on the torus. 
	\end{enumerate}
	\begin{proof} If at least one of $\Sigma/ \Delta^\pm$ is not transitive on  the corresponding set $\mathcal{G}^\pm$, then $\Sigma$ fixes a parallel class. 
		
		
		We now assume both $\Sigma/ \Delta^\pm$ are transitive on  $\mathcal{G}^\pm$. By Brouwer's Theorem,   $\Sigma/ \Delta^\pm$ is isomorphic and acts equivalently to either $\text{\normalfont SO}(2,\mathbb{R})$ or a finite covering group of $\text{\normalfont PSL}(2,\mathbb{R})$. Furthermore, it cannot be the case that both $\Sigma/ \Delta^\pm$ are isomorphic to $\text{\normalfont SO}(2,\mathbb{R})$, otherwise $\dim \Delta^\pm=2$, which in turn implies $\dim \Delta^+  \Delta^-  =4 >\dim \Sigma$, a contradiction. If $\Sigma/ \Delta^+$ is isomorphic to $\text{\normalfont PSL}^k(2,\mathbb{R})$, then    $\dim \Delta^ + =0$ and so $\Sigma$ is almost simple.  By Theorem \ref{almostsimplePSL},  $\Sigma$ is isomorphic to $\text{\normalfont \text{\normalfont PSL}}(2,\mathbb{R})$ and acts diagonally on the torus under a suitable coordinate system. The diagonal is a circle, as proved in Lemma \ref{actionalmostsimple}. Finally, $\Sigma$ cannot fix another circle because it only has one 1-dimensional orbit on the torus.\qedhere 
	\end{proof}
\end{lemma}

\begin{lemma} \label{3dfixedelements} Assume $\Sigma$ fixes at least one parallel class. Without loss of generality, let   $\pi$ be the fixed $(+)$-parallel class. Then exactly one of the following occurs.
	\begin{enumerate}
		\item $\Sigma$ fixes at least one point.

		\item  $\Sigma$ fixes no points but fixes and acts transitively on every (+)-parallel class. In this case $\Sigma \cong \text{\normalfont \text{\normalfont PSL}}(2,\mathbb{R})$ and $\mathbb{T}$ is isomorphic to a half-classical Minkowski plane $\mathcal{M}(f,id)$, where $f$ is an orientation-preserving homeomorphism of $\mathbb{S}^1$.

		\item  $\Sigma$ fixes no points but fixes and acts transitively on exactly one parallel class, which is $\pi$.  In this case  $\Sigma \cong  \text{\normalfont L}_2 \times \text{\normalfont SO}(2,\mathbb{R})$. 
		The factor group $\Sigma /\Delta^-$  
		is isomorphic   
		and  acts equivalently to $\text{\normalfont SO}(2,\mathbb{R})$
		on $\mathcal{G}^-$. 	Also,   $\Sigma /\Delta^+ $
		is isomorphic   
		and  acts equivalently to  $\text{\normalfont L}_2 $
		on   $\mathcal{G}^+\backslash\{\pi\}$.
		
	\end{enumerate}

	\begin{proof} Since $\Sigma$ fixes the $(+)$-parallel class $\pi$, $\Sigma / \Delta^+$ is not  transitive on  $\mathcal{G}^+$. If $\Sigma / \Delta^-$ is not  transitive on  $\mathcal{G}^-$, then  $\Sigma$ fixes at least one point.   
		
		We proceed by assuming  that $\Sigma/ \Delta^-$ is transitive on $\mathcal{G}^-$.  
		By Brouwer's Theorem,   $\Sigma/ \Delta^-$ is isomorphic and acts equivalently to either $\text{\normalfont SO}(2,\mathbb{R})$ or $\text{\normalfont PSL}^k(2,\mathbb{R})$ for some $k$. We consider these 2 cases separately.
		
		\begin{enumerate}[label=Case  \arabic*:,leftmargin=0pt,itemindent=*]

			\item $\Sigma/ \Delta^- \cong \text{\normalfont PSL}^k(2,\mathbb{R})$. Then $\dim \Delta^-=0$, and so $\Sigma$ is almost simple. From the first case in Theorem \ref{almostsimplePSL},  $\Sigma \cong \text{\normalfont \text{\normalfont PSL}}(2,\mathbb{R})$ and the plane $\mathbb{T}$ is determined. 
			
			\item $\Sigma/ \Delta^- \cong \text{\normalfont SO}(2,\mathbb{R})$.   
			If $\Sigma / \Delta^+$ acts trivially on $\mathcal{G}^+$, then $\Sigma = \Delta^+$.
			But since $\Delta^+ \cap \Delta^- = \{id\}$, it follows that $\Delta^- = \{id\}$, a contradiction to the dimensions of $\Sigma$ and $\Sigma/ \Delta^-$.  Hence $\Sigma / \Delta^+$ acts non-trivially on $\mathcal{G}^+$.  By Lemma \ref{salzmann9629}, there is an open subset $I \cong \mathbb{R}$ of $\mathcal{G}^+$ on which $\Sigma / \Delta^+$ acts transitively.  On the other hand,  $\Sigma$ is isomorphic to a subgroup of $\Sigma / \Delta^+ \times \Sigma / \Delta^-$, so that
			$$2 \le \dim \Sigma / \Delta^+ \le 3.$$ 
			From Brouwer's Theorem and Theorem \ref{almostsimplePSL}, it follows that $\Sigma / \Delta^+ \cong \text{\normalfont L}_2$.  
			By Lemma \ref{salzmann9312}, $\Sigma \cong \Sigma / \Delta^+ \times \Sigma / \Delta^- \cong  \text{\normalfont L}_2 \times \text{\normalfont SO}(2,\mathbb{R})$. 
			
			Since $\Delta^-$ fixes at most two $(+)$-parallel classes, so does $\Sigma$.
			Suppose   that	$\Sigma$ fixes exactly two $(+)$-parallel classes $\pi$ and  $\pi_2$.  Let $\pi_3$ be a $(+)$-parallel class different from $\pi$ and $\pi_2$. From the dimension formula, 
			$$
			2= \dim \Delta^- = \dim \Delta^-_{\pi,\pi_2} = \dim \Delta^-_{\pi,\pi_2,\pi_3} +\dim \Delta^-_{\pi,\pi_2}(\pi_3) \le 0+1=1, 
			$$
			a contradiction.

			Hence $\Sigma$ fixes exactly one $(+)$-parallel class $\pi$. In particular, $I=\mathcal{G}^+\backslash\{\pi\}$ and $\Sigma / \Delta^+$ acts equivalently to $\text{\normalfont L}_2$ on $\mathcal{G}^+\backslash\{\pi\}$, by Brouwer's Theorem.  
			\qedhere  
		\end{enumerate} 
	\end{proof}
	
\end{lemma}

\begin{lemma} \label{3dfixepoint} Assume  $\Sigma$ fixes at least one point.  Under a suitable coordinate system, let  $p=(\infty,\infty)$ be a fixed point. Then exactly one of the following occurs.
	\begin{enumerate}
		\item $\Sigma$ fixes exactly two parallel points. In this case $\Sigma \cong \Phi_d$, for some $d \le 0$. The coordinates may be chosen such that the second fixed point is $(0,\infty)$ and the  action of $\Sigma$ is described by the maps
		$$
		\{ (x,y) \mapsto (ax,by+c) \mid a,b>0,c \in \mathbb{R} \},
		$$
		when $\Sigma \cong \Phi_0$, and  
		$$
		\{ (x,y) \mapsto (a \text{ } \text{sgn}(x) \cdot |x|^b ,b^dy+c) \mid a,b>0,c \in \mathbb{R} \},
		$$
		when $\Sigma \cong \Phi_d$, for some  $d<0$. 
		
		\item$\Sigma$ fixes exactly one point, which is $p$.  In this case the derived plane $\mathbb{T}_p$ is Desarguesian and $\Sigma \cong \Phi_d$, for some $d \in \mathbb{R} \cup \{ \infty\}$. The coordinates may be chosen such that the action of  $\Sigma$ is the standard action of $\Phi_d$ on $\mathbb{R}^2$.
	\end{enumerate}	
	\begin{proof} There are 3 cases depending on the transitivity of $\Sigma/ \Delta^\pm$ on $\mathcal{G}_\pm\backslash \{ [p]_\pm\}$.
		\begin{enumerate}[label=Case  \arabic*:,leftmargin=0pt,itemindent=*]
			\item Neither of  $\Sigma/ \Delta^\pm$ is transitive on   $\mathcal{G}_\pm \backslash \{ [p]_\pm\}$. This implies there is an additional  fixed point $q$ nonparallel to $p$.  From the dimension formula, for a point $r$ nonparallel to $p$ and $q$, we have
			\begin{align*}
			3=\dim \Sigma_p 	= \dim \Sigma_{p,q}
			&=  \dim \Sigma_{p,q,r}+\dim \Sigma_{p,q}(r)  \le 2, 
			\end{align*} 
			a contradiction. Hence this case cannot occur.

			\item  $\Sigma/ \Delta^-$ is transitive on $\mathcal{G}^- \backslash \{ [p]_- \}$, $\Sigma/ \Delta^+$ is not transitive on $\mathcal{G}^+ \backslash \{ [p]_+\}$, or vice versa.  By Lemma \ref{salzmann9629}, $\Sigma$ fixes  an additional point $q \in [p]_-$. We  show   $\Sigma$ fixes at most two  points. By changing the coordinates if necessary,  we can assume $q=(0,\infty)$.  Suppose to the contrary that $\Sigma$  fixes three   points on $[p]_-$.
			Then  $\Delta^-$ fixes  three (+)-parallel classes pointwise  and so must be trivial. By Brouwer's Theorem,   $\Sigma \cong \Sigma/ \Delta^- \cong \widetilde{\text{\normalfont \text{\normalfont PSL}}(2,\mathbb{R})}$, which contradicts Theorem \ref{almostsimplePSL}.
			
			Hence $\Sigma$  fixes precisely two parallel points $p$ and $q$.  
			On the derived plane $\mathbb{T}_p$, $\Sigma$ induces a group of automorphisms that  fixes precisely a line. By \cite[ Theorem 7.5B]{groh1983},  $\Sigma$ is isomorphic to either  $\mathbb{R} \times \text{\normalfont L}_2 \cong \Phi_0$ or $\Phi_d$, for some $d<0$, and the action of the group is described as in the statement of the lemma. 
			
			\item  Both  $\Sigma/ \Delta^\pm$ are   transitive on   $\mathcal{G}^\pm\backslash \{ [p]_\pm\}$. In this case, $p$ is the only point fixed by $\Sigma$.  Brouwer's Theorem implies that $\Sigma / \Delta^\pm$ is isomorphic to $\mathbb{R}$, $\text{\normalfont L}_2$ or $\widetilde{\text{\normalfont \text{\normalfont PSL}}(2,\mathbb{R})}$.    From Theorem \ref{almostsimplePSL}, we can rule out the last case.    Since $\Sigma$ is isomorphic to a subgroup of $\Sigma / \Delta^+ \times \Sigma / \Delta^-$,    it cannot be the case that both $\Sigma / \Delta^\pm \cong \mathbb{R}$. This leads to 2 subcases. 
			
			\begin{enumerate}[label=Subcase  3\Alph*:,leftmargin=0pt,itemindent=*]
				\item  $\Sigma / \Delta^+ \cong \mathbb{R} $ and $\Sigma / \Delta^- \cong \text{\normalfont L}_2$, or vice versa. 
				By Lemma \ref{salzmann9312}, $\Sigma = \Sigma / \Delta^+ \times \Sigma / \Delta^- \cong \mathbb{R} \times \text{\normalfont L}_2$.  By Brouwer's Theorem, the action of $\Sigma / \Delta^\pm$ is standard, and so the action of  $\Sigma$ on $\mathcal{P}$ is described by the maps
				$$ \{ (x,y) \mapsto (x+b,ay+c)\mid a>0,b,c \in \mathbb{R} \},$$ 
				in suitable coordinates. In particular, $\Sigma \cong \Phi_\infty$. 	
				
				When 		 $\Sigma / \Delta^- \cong \mathbb{R}  $ and $\Sigma / \Delta^+ \cong \text{\normalfont L}_2 $, we obtain $\Sigma \cong \Phi_0$ in a similar fashion. 	
				
				%
				%

				\item   $\Sigma / \Delta^\pm \cong \text{\normalfont L}_2$.  We show that $\Sigma \cong \Phi_d$ for some $d \in \mathbb{R} \cup \{ \infty\}$.   
				
				Let  $\overline{\mathcal{P}}:=\mathcal{P} \backslash ([p]_+ \cup [p]_-)  \cong \mathbb{R}^2$.  We consider the action of $\Delta^+\Delta^-$ on $\overline{\mathcal{P}}$. We have $\dim \Delta^+ = 1$. If $\Delta^+$ fixes a parallel class $[q]_- \in \mathcal{G}^-\backslash \{ [p]_-\}$, then, as a normal subgroup, it fixes the orbit of $[q]_-$ pointwise. But since  $\Sigma/ \Delta^-$ is transitive on $\mathcal{G}^-\backslash \{ [p]_-\}$, $\Delta^+$ is then trivial, which is impossible. Hence $\Delta^+$ is transitive on $\mathcal{G}^-\backslash \{ [p]_-\}$. By Brouwer's Theorem, $\Delta^+$
				is isomorphic and acts  equivalently to  $\mathbb{R}$. 
				With the same reasoning for $\Delta^-$, it follows  that $\Delta^+\Delta^-$ is the translation group $\mathbb{R}^2$ (in suitable coordinates) and is sharply transitive on $\overline{\mathcal{P}}$. 
				
				Denote $o \coloneqq(0,0)$. By \cite[Proposition 91.2]{salzmann1995},  $\Sigma=\Delta^+\Delta^-\Sigma_o$.     Since $\Delta^+\Delta^-$ is a normal subgroup of $\Sigma$ and $\Delta^+\Delta^- \cap \Sigma_o$ is trivial,   $\Sigma=\Sigma_o \ltimes \Delta^+\Delta^-$.  
				From  the dimensions of $\Sigma$ and $\Delta^\pm$,   $\dim \Sigma_o =1$. By \cite[Corollary 94.39]{salzmann1995}, $\Sigma_o \cong \mathbb{R}$. The action of  $\Sigma_o$ is then described by the maps
				$$ \{ (x,y) \mapsto (x,ay) \mid a>0   \},  $$ or
				$$ \{ (x,y) \mapsto (ax,a^dy)\mid a>0   \}.$$   
				This shows that $\Sigma \cong \Phi_d$ for some $d \in \mathbb{R} \cup \{ \infty\}$.   
				
				In both subcases,  $\Delta^+ \times \Delta^-$ contains a normal subgroup that is sharply transitive on $\overline{\mathcal{P}}$. 
				Hence $\Sigma$ induces a 3-dimensional group of automorphisms  acting transitively on the point set  of  $\mathbb{T}_p$.  The list of possibilities for $\mathbb{T}_p$ is given by \cite[Main Theorem 2.6]{groh1982a} and a case by case check shows that $\mathbb{T}_p$ is Desarguesian, cf. \cite[Subsection 7.1.1]{duythesis}. 
				\qedhere
				
			\end{enumerate}

		\end{enumerate}
	\end{proof}
\end{lemma}

 \pagebreak

\bibliography{bib2} 
\bibliographystyle{abbrv}


\end{document}